\documentclass{article}
\usepackage[utf8]{inputenc}
\usepackage[english]{babel}

\usepackage{csquotes, soul}
\usepackage{amsmath}
\usepackage{amsthm}
\usepackage{amssymb}
\usepackage{amsfonts}
\usepackage{amscd}
\usepackage{hyperref}
\usepackage{enumerate}
\usepackage{fancyhdr}
\usepackage{mathrsfs}
\usepackage{polynom}
\usepackage{graphicx}
\usepackage{tikz}
\usepackage{hyperref}
\usepackage{verbatim}
\usepackage{accents}
\usepackage{comment}

\hypersetup{breaklinks=true}
 
\usepackage{geometry}
\newcommand{\arc}{\mathcal{A}(S)} 
\newcommand{\arcp}{\mathcal{A}(S')}

\newcommand{\flip}{\mathcal{F}(S)}
\newcommand{\flipp}{\mathcal{F}(S')}

\newcommand{\emod}{\text{Mod}^\pm(S)}

\theoremstyle{plain}
\newtheorem{theorem}{Theorem}[section]

\newtheorem{corollary}[theorem]{Corollary}

\theoremstyle{remark}

\makeatletter
\def\blfootnote{\xdef\@thefnmark{}\@footnotetext}
\makeatother

\begin{document}

\title{Flip graph and arc complex finite rigidity}
%\author{Chandrika Sadanand and Emily Shinkle}
%\author[1]{Chandrika Sadanand%\corref{cor1}
%}
%\ead{c.sadanand@bowdoin.edu}
%\author[2]{Emily Shinkle %\fnref{fn2}
%}
%\ead{esshinkle@gmail.com}
%\cortext[cor1]{Corresponding author}
%\affiliation[1]{organization={Bowdoin College},
%addressline={255 Maine St.},
%postcode={04011},
%city={Brunswick, ME},
%country={USA}}
%\affiliation[2]{organization={University of Illinois Urbana-Champaign},
%addressline={1409 W. Green St.},
%postcode={61801},
%city={Urbana, IL},
%country={USA}}

%\fntext[fn2]{Emily's permanent address}

\date{}
%\begin{keyword}
%\MSC[2020]{57K20, 57M50, 05C25, 20F65} \sep %Mapping class groups \sep Arc complex \sep Rigidity
%\end{keyword}
\author{Chandrika Sadanand and Emily Shinkle}
\maketitle
\begin{abstract} 
A subcomplex $\mathcal{X}$ of a cell complex $\mathcal{C}$ is called \emph{rigid} with respect to another cell complex $\mathcal{C}'$ if every injective simplicial map $\lambda:\mathcal{X} \rightarrow \mathcal{C}'$ has a unique extension to an injective simplicial map $\phi:\mathcal{C}\rightarrow \mathcal{C}'$. We say that a cell complex exhibits \emph{finite rigidity} if it contains a finite rigid subcomplex. Given a surface with marked points, its \textit{flip graph} and \textit{arc complex} are simplicial complexes indexing the triangulations and the arcs between marked points, respectively. In this paper, we leverage the fact that the flip graph can be embedded in the arc complex as its dual to show that finite rigidity of the flip graph implies finite rigidity of the arc complex. Thus, a recent result of the second author on the finite rigidity of the flip graph implies finite rigidity of the arc complex for a broad class of surfaces. Notably, this includes surfaces with boundary -- a setting where finite rigidity of the arc complex was previously unknown.
\end{abstract}

%\maketitle 
%\blfootnote{MSC2020: 57K20, 57M50, 05C25, 20F65}

\section{Introduction}\label{section: introduction}

The curve graph $\mathcal{C}(S)$ of a surface $S$ is a graph that indexes the simple closed curves on $S$. Its study has provided great insights into the geometry of the Teichm\"{u}ller space and mapping class group of $S$ (e.g.~\cite{harer1986virtual}, \cite{Harvey}). Many of these insights are possible because the induced homomorphism $f:\text{Mod}^{\pm}(S) \to \text{Aut}(\mathcal{C}(S))$ is an isomorphism (except in a few well-understood cases). This property of $\mathcal{C}(S)$ is known as \emph{simplicial rigidity} or \emph{rigidity}. It was proved by Ivanov~\cite{ivanov1997automorphisms} for surfaces of genus two or more, without or with boundary. Korkmaz \cite{korkmaz1999automorphisms} and Luo \cite{luo1999automorphisms} settle the cases of genus zero and one.

Aramayona and Leininger were able to strengthen these results considerably by showing that (except for the twice punctured torus) there exists a finite subcomplex $\mathcal{X}$ of $\mathcal{C}(S)$ such that every locally injective map $\mathcal{X} \to \mathcal{C}(S)$ is the restriction of $f(g)$ to $\mathcal{X}$ for some (unique in most cases) $g \in \text{Mod}^{\pm}(S)$~\cite{aramayona2013finite}.
%They further show that for every orientable surface, there is a finite subcomplex $\mathcal{X}$ of $\mathcal{C}(S)$, so that  every locally injective map $\lambda:\mathcal{X} \rightarrow \mathcal{C}(S)$ has a unique extension to an automorphism of $\mathcal{C}(S)$.
This is referred to as \emph{finite rigidity}. Similar rigidity results have since been proved for many other related complexes and used to provide new proofs of existing results (e.g.~\cite{disarlo2015combinatorial}) as well as to prove new theorems (e.g.~\cite{disarlo2019geometry}).

Of interest in this paper are the arc complex $\arc$ and the flip graph $\flip$. They respectively index the arcs and triangulations of a surface $S$ with marked points. Both complexes were introduced by Harer~\cite{harer1986virtual} to study homological properties of the mapping class group. Using distinct methods, the second author has recently proved the finite rigidity of $\flip$ for surfaces $S$ without or with boundary \cite{shinkle2022finite} 
%to show finite rigidity of $\arc$ (Corollary \ref{cor: main result}).
as well as finite rigidity of $\arc$
%was previously proved by the second author
for surfaces $S$ without boundary only \cite{shinkle2020finite}.

We exploit the duality of $\arc$ and $\flip$ to show that finite rigidity of $\flip$ implies finite rigidity of $\arc$ (Theorem \ref{thm: new result}). Thus, in addition to demonstrating the utility of the flip graph--arc complex duality and providing a novel method for proving the finite rigidity of $\arc$ for surfaces without boundary, we show that finite rigidity of $\arc$ also holds for surfaces with boundary (Corollary \ref{cor: main result}). This was previously unknown.

\paragraph{Acknowledgements} We would like to thank Yair Minsky for suggesting the exploration of duality of the arc complex and the flip graph, and its implications on finite rigidity.

\section{Background}\label{section: background}

Let $S$, $S'$ be compact, connected, orientable surfaces, possibly with boundary. Assume $S$ and $S'$ each contain a fixed, finite, nonzero number of marked points, with at least one marked point on each boundary component. We write $S_{g,n,(p_1,\ldots,p_b)}$ for a surface of genus $g$ with $n$ marked points in its interior and $b$ boundary components where the number of marked points on the $i^{th}$ boundary component is $p_i$. If $b  =  0$, we write $S_{g,n}$. Let $d(S) = 6g + 3n + 3b + p_1 + p_2 + \cdots + p_b - 6$, which is the number of arcs needed to triangulate $S$ in most cases. The \emph{extended mapping class group} of $S$, $\emod$, is the group of homeomorphisms of $S$ up to isotopy.

The arc complex $\arc$ of a surface $S$ with marked points is the simplicial complex whose 0-simplices represent isotopy classes of essential arcs on $S$ and whose $k$ simplices join collections of $k+1$ vertices representing pairwise disjoint classes of arcs for $k\geq 1$. Further details and examples are available in %\cite{placeholder}.
\cite{shinkle2020finite}.

Arc complex rigidity was proved for surfaces without boundary (Irmak-McCarthy~\cite{irmak2010injective}) and surfaces with boundary (Disarlo~\cite{disarlo2015combinatorial}). In 
\cite{shinkle2020finite},
%\cite{placeholder},
the second author improved this result by proving finite rigidity. 

\begin{theorem}[Arc complex finite rigidity]
\label{thm: arc complex finite rigidity}
Let $S=S_{g,n}$ where $S\not \cong S_{0,n}$ for $n\leq 3$ or $S_{1,1}$. There exists a finite subcomplex $\mathcal{X}\subseteq \arc$, such that for any $S'$ with $d(S) \geq d(S')$, and for any locally injective simplicial map $\lambda:\mathcal{X}\rightarrow \arcp$, there exists a 
homeomorphism $h:S\rightarrow S'$ which induces $\lambda$, unique up to isotopy.
%simplicial isomorphism $\phi:\arc \rightarrow \arcp$ such that $\lambda = \phi|_{\mathcal{X}}$. 
\end{theorem}

The flip graph $\flip$ of $S$ is the graph whose vertices represent triangulations of $S$, and whose edges connect pairs of vertices whose corresponding triangulations differ by a single arc. More details and examples are available in %\cite{placeholder}.
\cite{shinkle2022finite}.

Flip graph rigidity results exist for surfaces without boundary (Korkmaz--Papadopoulis~\cite{korkmaz2012ideal}) and with boundary (Aramayona--Koberda--Parlier~\cite{aramayona2015injective}).
Again, the second author improved this result by proving finite rigidity in 
\cite{shinkle2022finite}.
%\cite{placeholder}.

\begin{theorem}[Flip graph finite rigidity]
\label{thm: flip graph finite rigidity}
Let $S=S_{g,n,(p_1,\ldots,p_b)}$ where 
$S$ cannot be embedded in $S_{0,n}$ for $n\leq 4$ and $S \not \cong S_{1,n}$ with $n\leq 2$.
% $S\not \cong $ $S_{0,0,(1,1)}, S_{1,1}$. 
Let $S'$ be any surface with $d(S)=d(S')$. Then there exists a finite subgraph $\mathcal{X}\subseteq \flip$ such that any injective homomorphism $\lambda:\mathcal{X} \rightarrow \flipp$ is induced by a homeomorphism $h:S \rightarrow S'$, unique up to isotopy.
\end{theorem}

% Such an $\mathcal{X}$ is called a \textit{isomorphism-inducing finite rigid subcomplex} of $\arc$. Due to \cite{irmak2010injective}, if $S\not \cong S_{0,3}$, then the existence of such an isomorphism is equivalent to the existence of a homeomorphism $h:S\rightarrow S'$ where $h(a)=\lambda(a)$ for any arc $a$ represented by a vertex in $\mathcal{X}$.

% Such an $\mathcal{X}$ is called a \textit{finite rigid subgraph} of $\flip$. Again, this injective simplicial map can be realized by an embedding (not necessarily unique) $i:S\rightarrow S'$ and by a homeomorphism (unique up to isotopy) $h:S\rightarrow S'$ if $d(S) = d(S')$ (see \cite{korkmaz2012ideal}, \cite{aramayona2015injective}).

% Observe that in Theorem \ref{thm: arc complex finite rigidity}, the subcomplex $\mathcal{X}$ depends only on $S$, whereas in Theorem \ref{thm: flip graph finite rigidity}, the subcomplex $\mathcal{X}$ depends both on $S$ and $S'$. If $S \cong S_{0,0,(1,1)}$, then Theorem \ref{thm: flip graph finite rigidity} holds for some $S'$ but not others. For $S\cong S_{1,1}$, the result is not known. See 
% %\cite{shinkle2022finite} 
% \cite{placeholder} for more details.
\section{Duality and finite $\mathcal{A}(S)$ rigidity for surfaces with boundary}
Every vertex $T$ in $\flip$ corresponds to a triangulation of $S$, which corresponds to a top-dimensional simplex $\tau$ of $\arc$. Notice that if two top-dimensional simplices $\tau, \tau'$ of $\arc$ are adjacent, then they intersect in a codimension one simplex. This implies that they differ by a single arc, in other words, by a flip. Thus the vertices of $\flip$ that $\tau$ and $\tau'$ correspond to are adjacent in $\flip$. This allows $\flip$ to be embedded in $\arc$ as its dual graph. Define the bijection $\pi:\flip \rightarrow \{\text{top-dimensional simplices of } \arc\}$ which sends a triangulation $T$ to the simplex in $\arc$ whose vertices correspond to the arcs which compose $T$.

In this note, we prove that a variation of Theorem~\ref{thm: arc complex finite rigidity} is implied by Theorem~\ref{thm: flip graph finite rigidity}. In other words, with some hypotheses, the existence of a finite rigid subgraph of the flip graph implies the existence of a finite rigid subcomplex of the arc complex.

% \Restate{Theorem}{thm: new result}{}{Suppose $d(S) = d(S') \geq 2$ and that $\flip$ is not finite. Further suppose $\mathcal{X}_{\mathcal{F}}\subseteq \flip$ is a finite rigid subgraph containing %contains
% a 2-ball in $\flip$.
% % and that for any injective simplicial map $\lambda_{\mathcal{F}}:\mathcal{X}_{\mathcal{F}}\rightarrow \flipp$, there exists a unique injective simplicial map $\phi_{\mathcal{F}}:\flip \rightarrow \flipp$ such that $\lambda_{\mathcal{F}} = \phi_{\mathcal{F}}|_{\mathcal{X}_{\mathcal{F}}}$.
% % Let 
% % $\mathcal{X}_{\mathcal{A}} = \cup_{T\in \mathcal{X}_{\mathcal{F}}} \pi(T)$. 
% Then $\mathcal{X}_{\mathcal{A}} = \cup_{T\in \mathcal{X}_{\mathcal{F}}} \pi(T)$ is a finite rigid subcomplex of $\arc$.
% %for any injective simplicial map $\lambda_{\mathcal{A}}:\mathcal{X}_{\mathcal{A}} \rightarrow \arcp$, there exists a unique simplical injective simplicial map $\phi_{\mathcal{A}}:\arc \rightarrow \arcp$ such that $\phi_{\mathcal{A}}|_{\mathcal{X}_{\mathcal{A}}} = \lambda_{\mathcal{A}}$.
% }
By $n$-ball, we refer to the collection of vertices in a graph (or subgraph containing exactly those vertices and their connecting edges) which are distance at most $n$ from a fixed center vertex for the ball. An $n$-shell refers to only those points which are distance exactly $n$ from the ball center. Note that $n$ does not indicate any sort of dimension in this context.

\begin{theorem}[Flip graph finite rigidity implies arc complex finite rigidity]\label{thm: new result}
Suppose $d(S) = d(S') \geq 2$ and that $\flip$ is not finite. Further suppose $\mathcal{X}_{\mathcal{F}}\subseteq \flip$ contains a 2-ball in $\flip$ and that for any injective simplicial map $\lambda_{\mathcal{F}}:\mathcal{X}_{\mathcal{F}}\rightarrow \flipp$, there exists a unique (up to isotopy) homeomorphism 
$h_{\mathcal{F}}:S \rightarrow S'$ which induces $\lambda_{\mathcal{F}}$.
%injective simplicial map $\phi_{\mathcal{F}}:\flip \rightarrow \flipp$ such that $\lambda_{\mathcal{F}} = \phi_{\mathcal{F}}|_{\mathcal{X}_{\mathcal{F}}}$. 
Let 
$\mathcal{X}_{\mathcal{A}} = \cup_{T\in \mathcal{X}_{\mathcal{F}}} \pi(T)$. Then for any injective simplicial map $\lambda_{\mathcal{A}}:\mathcal{X}_{\mathcal{A}} \rightarrow \arcp$, 
there exists a unique (up to isotopy) homeomorphism $h_{\mathcal{A}}:S \rightarrow S'$ which induces $\lambda_{\mathcal{A}}$.
%there exists a unique simplical injective simplicial map $\phi_{\mathcal{A}}:\arc \rightarrow \arcp$ such that $\phi_{\mathcal{A}}|_{\mathcal{X}_{\mathcal{A}}} = \lambda_{\mathcal{A}}$.
\end{theorem}

Note that if $\mathcal{X}_{\mathcal{F}}$ is finite, then $\mathcal{X}_{\mathcal{A}}$ is as well. Theorem \ref{thm: new result} together with Theorem \ref{thm: flip graph finite rigidity} implies the main result of this paper --- that finite rigidity holds for the arc complexes of all but a few exceptional surfaces. The exclusions are those surfaces excluded in Theorem \ref{thm: flip graph finite rigidity} or Theorem \ref{thm: new result}. For the former, we exclude surfaces $S$ which can be embedded in $S_{0,n}$ for $n\leq 4$ as well as $S_{1,n}$ for $n\leq 2$. For the later, in all cases where $d(S) < 2$, the surface $S$ can be embedded in $S_{0,n}$ for $n\leq 4$, and $S$ is already excluded. Finally, we must exclude any remaining surfaces $S$ for which $\flip$ is finite: $S_{0,0,(p_1)}$ for $p_1 > 4$. We will refer to these surfaces as \emph{exceptional} surfaces.

\begin{corollary}[Arc complex finite rigidity extended]

\label{cor: main result}

Let $S=S_{g,n,(p_1,\ldots,p_b)}$ be a surface that is not exceptional and let $d(S)=d(S')$. Then there exists a finite subcomplex $\mathcal{X}\subseteq \arc$, such that for any injective simplicial map $\lambda:\mathcal{X}\rightarrow \arcp$, there exists a unique (up to isotopy) homeomorphism $h_{\mathcal{A}}:S \rightarrow S'$ which induces $\lambda_{\mathcal{A}}$.
%there exists a unique injective simplicial map $\phi:\arc \rightarrow \arcp$ such that $\lambda = \phi|_{\mathcal{X}}$.

\end{corollary}

Notably, this variation on Theorem \ref{thm: arc complex finite rigidity} applies to surfaces with boundary, where finite rigidity of the arc complex was previously unknown. Further, the techniques through which the new result is obtained are substantially different than those used for the original theorem. 
% Note that Theorem \ref{thm: flip graph finite rigidity} applies to surfaces with boundary, whereas Theorem \ref{thm: arc complex finite rigidity} does not. So we can conclude that Theorem \ref{thm: arc complex finite rigidity} holds for a larger class of surfaces than previously known, albeit in a slightly weaker form. The original proof of Theorem \ref{thm: arc complex finite rigidity} was constructive, i.e. the subcomplex $\mathcal{X}$ is explicitly defined. In contrast, the proof of Theorem \ref{thm: flip graph finite rigidity} is not constructive, so we are not provided with an exact description of the finite rigid subcomplexes in $\arc$ for surfaces with boundary.

% \section{A variation on Shinkle's arc complex finite rigidity}
\section{Proofs}

\begin{proof}[Proof of Theorem~\ref{thm: new result}]
Let $\mathcal{X}_{\mathcal{A}}$ be defined as in Theorem~\ref{thm: new result}. Suppose $\lambda_{\mathcal{A}}:\mathcal{X}_{\mathcal{A}} \rightarrow \arcp$ is an injective simplicial map. Define a map $\lambda_{\mathcal{F}}:\mathcal{X}_{\mathcal{F}}\rightarrow \flipp$ which sends a triangulation $T \in \mathcal{X}_{\mathcal{F}}$ to $\cup_{\{\text{arcs } a \in T\}} \lambda_{\mathcal{A}}(a)$, which is a triangulation of $S'$ (since $d(S) =d(S')$). Recall that $\lambda_{\mathcal{A}}$ is indeed defined for all $a\in T$ due to the construction of $ \mathcal{X}_{\mathcal{A}}$. We can see that $\lambda_{\mathcal{F}}$ is injective because $\lambda_{\mathcal{A}}$ is. Suppose $T, T' \in \mathcal{X}_{\mathcal{F}}$ are adjacent. This can be characterized by the fact that $T\backslash T' = \{a_0\}$ and $T'\backslash T = \{a_1\}$ for some arcs $a_0, a_1$ on $S$. Then $\lambda_{\mathcal{F}}(T)\backslash \lambda_{\mathcal{F}}(T') = \{\lambda_{\mathcal{A}}(a_0)\}$ and $\lambda_{\mathcal{F}}(T')\backslash \lambda_{\mathcal{F}}(T) = \{\lambda_{\mathcal{A}}(a_1)\}$, so $\lambda_\mathcal{F}(T)$ and $\lambda_{\mathcal{F}}(T')$ are also adjacent. This implies that $\lambda_{\mathcal{F}}$ is simplicial.

By assumption, 
%there is a unique injective simplicial map $\phi_{\mathcal{F}}:\flip \rightarrow \flipp$ such that $\lambda_{\mathcal{F}} = \phi_{\mathcal{F}}|_{\mathcal{X}_{\mathcal{F}}}$. By \cite{korkmaz2012ideal}, there is a unique homeomorphism $h:S\rightarrow S'$ which agrees with $\phi$ on every triangulation $T$ of $S$, and thus $h(T)=\lambda_{\mathcal{F}}(T)$ for all $T\in \mathcal{X}_{\mathcal{F}}$. (Note that this means $S\cong S'$.) Let $\phi_{\mathcal{A}}:\arc \rightarrow \arcp$ be the map induced (uniquely) by $h$. Then it suffices to show that $\lambda_{\mathcal{A}}(a) = \phi_{\mathcal{A}}(a)$ for all $a\in \mathcal{X}_{\mathcal{A}}$.
there is a unique (up to isotopy) homeomorphism $h:S\rightarrow S'$ which induces $\lambda_{\mathcal{F}}$. More precisely, for any $T \in \mathcal{X}_{\mathcal{F}}$, $\lambda_{\mathcal{F}}(T) = h(T)$. We claim that $h$ also uniquely induces $\lambda_{\mathcal{A}}$. To show that $h$ induces $\lambda_{\mathcal{A}}$, it suffices to prove that for all $a\in\mathcal{X}_{\mathcal{A}}$, that $\lambda_{\mathcal{A}}(a)=h(a)$.

As we will see, there is no arc which is contained in every vertex of a 2-ball in $\flip$. By assumption, we are only considering the cases where $\flip$ is infinite (the generic case). If $\flip$ is infinite then it is also infinite diameter, since there is a finite upper bound on the degrees of its vertices. Thus there are indeed points in the 2-shell of this 2-ball. Call the center of the 2-ball $T$. Any arc $a$ which is in every element on the 1-shell of $T$ must be unflippable in $T$. This means that $a$ is the inner arc of a folded triangle in $T$. But then if $b$ is the outer arc of this folded triangle, then the triangulation resulting from a flip of $T$ along $b$ has $a$ as a flippable arc, and the subsequent flip of this triangulation along $a$ is in the 2-shell of $T$ and does not contain $a$. See 
\cite{shinkle2022finite} 
%\cite{placeholder}
for a more thorough explanation of this terminology.

Since $\mathcal{X}_{\mathcal{F}}$ contains a 2-ball and is connected by hypothesis, given an arc $a$, we may find adjacent triangulations $T, T'\in \mathcal{X}_{\mathcal{F}}$ such that $a \in T$ and $a \notin T'$. Then we see that 
\[
%\phi_{\mathcal{A}}(a) = 
\lambda_{\mathcal{A}}(a) = \lambda_{\mathcal{F}}(T)\backslash \lambda_{\mathcal{F}}(T') = h(T)\backslash h(T') = h(a).\]

%Suppose another injective simplicial map $\phi_{\mathcal{A}}'$ agrees with $\lambda_{\mathcal{A}}$ on $\mathcal{X}_{\mathcal{A}}$. If $S\cong S_{0,3}$, then $\flip$ has diameter 2 and hence $\mathcal{X}_{\mathcal{F}} = \flip$. This implies that $\mathcal{X}_{\mathcal{A}} = \arc$, hence $\phi_{\mathcal{A}}' = \phi_{\mathcal{A}}$. Otherwise, consider that $\phi_{\mathcal{A}}'$ is induced by a homeomorphism $h':S\rightarrow S'$ by \cite{irmak2010injective}. Then since $h$ and $h'$ agree on each arc in a triangulation of $S$, they are isotopic and hence $\phi_{\mathcal{A}} = \phi_{\mathcal{A}}'$. In either case, the extension is unique.
Finally, $h$ is unique because any two homeomorphisms which agree on all of $\lambda_{\mathcal{A}}$ agree on a set of arcs which form a triangulation and thus the homeomorphisms are isotopic.

\end{proof}

\begin{proof}[Proof of Corollary~\ref{cor: main result}]

Let $S=S_{g,n,(p_1,\ldots,p_b)}$ not exceptional and $d(S)=d(S')$. Let $\mathcal{X}_{\mathcal{F}}$ be the finite subgraph from Theorem \ref{thm: flip graph finite rigidity}. If $\mathcal{X}_{\mathcal{F}}$ does not contain a 2-ball in $\flip$, then using the technique employed in the proof of Corollary 1.4 in~%\cite{placeholder},
\cite{shinkle2022finite},
we can expand $\mathcal{X}_{\mathcal{F}}$ so that it does contain a 2-ball while retaining the finite rigidity property. Finally, let $\mathcal{X}_{\mathcal{A}}$ be the corresponding finite subcomplex of $\arc$ from Theorem~\ref{thm: new result}. This subcomplex $\mathcal{X}_{\mathcal{A}}$ satisfies the conclusion of the Corollary in question.

\end{proof}

\section{Discussion}

Interestingly, the converse of Theorem~\ref{thm: new result} does not hold --- finite rigidity of the arc complex does not imply finite rigidity of the flip graph. To understand why, consider the case of $S_{1,1}$, the torus with one marked point. Theorem \ref{thm: arc complex finite rigidity} holds for $S_{1,1}$, whereas Theorem \ref{thm: flip graph finite rigidity} does not, because the uniqueness of the extension is not guaranteed. The arc complex is the well-known Farey graph, which does exhibit finite rigid subcomplexes; whereas its flip graph is a trivalent tree, which does not. The proof of our result above relies on the ability to induce a map on some portion of $\flip$ from a map on some portion of $\arc$. This is simple because if one knows the image of each arc in a triangulation, they know the image of the triangulation. But the converse is not true. The mapping of one triangulation to another provides no information about a mapping between the individual arcs which compose each triangulation. This is not to say that such a task cannot be done provided sufficient additional information. The results of \cite{korkmaz2012ideal} and \cite{aramayona2015injective} both employ such a technique to induce maps defined on all of $\arc$ from maps defined on all of $\flip$. However, this cannot be done in general for any map defined on some portion of $\flip$.

\bibliographystyle{plain}
\bibliography{main}
%\printbibliography

\end{document}